\documentclass[11pt]{amsart}
\usepackage{amssymb, amsmath}
\usepackage{amsthm}
\usepackage{amscd}
\usepackage{graphicx}
\newtheorem{theorem}{Theorem}

\newtheorem{question}[theorem]{Question}

\newtheorem{conjecture}[theorem]{Conjecture}
\theoremstyle{definition}

\newtheorem*{remark}{Remark}
\newtheorem*{acknowledgement}{Acknowledgement}

\title[On $\varepsilon$-maps in dimension 1]
{A remark on $\varepsilon$-maps in dimension 1}
\author{T. T$\hat{\mathrm{a}}$m Nguy$\tilde{\hat{\mathrm{e}}}$n-Phan}
\address{Max Planck Institute for Mathematics\\
Bonn\\
Germany}
\email{tam@mpim-bonn.mpg.de}

\def\beqa{\begin{eqnarray}}
\def\eeqa{\end{eqnarray}}
\def\beqa{\begin{eqnarray}}
\def\eeqa{\end{eqnarray}}

\input xy
\xyoption{all}

\begin{document}

\begin{abstract}
Let $f\colon \mathbb{S}^1\rightarrow G$ be a surjective map from the standard unit circle to a graph $G$ such that the pre-image of each point has diameter less than $\varepsilon$. If $\varepsilon$ is small enough, does $f$ split as a free factor in $\pi_1(G)$?
\end{abstract}
\maketitle

\section{Introduction}
Bestvina-Brady's paper ``Morse theory and finiteness propeties of groups" ended with the following conjecture (\cite{bestvinabrady}).

\begin{conjecture}[Bestvina-Brady]\label{conjecture}
Let $L$ be a finite 2-complex. Fix a metric on $L$. There is $\varepsilon>0$ such that if $g\colon L\rightarrow K$ is a surjective, PL $\varepsilon$-map\footnote{A map is an $\varepsilon$-map if pre-images of points have diameters less than $\varepsilon$.}, then $K$ is homotopy equivalent to $L$ with $1$-cells and $2$-cells attached.
\end{conjecture}

They had a specific complex $L$ in the paper and it is a spine of the Poincare homology sphere, and they remarked that if the conjecture is true for this 2-complex, then the Eilenberg-Ganea conjecture is false. In this paper we will take the point of view that Conjecture \ref{conjecture} is clearly more interesting.

While Conjecture \ref{conjecture} might sound ``obviously true", it remains open. If we reduce the dimension of $L$ from $2$ to $1$, then Conjecture \ref{conjecture}, in its literal form, is true for a rather dull reason, which is that such $PL$-quotients of a circle by $\varepsilon$-maps are graphs that are not trees, and such graphs are always circles with 1-cells attached. 

A more interesting version of Conjecture \ref{conjecture} in dimension $1$ is the following.

\begin{question}\label{question}
Is there $\varepsilon >0$ small enough such that if $f\colon \mathbb{S}^1 \rightarrow G$ is a surjective, PL $\varepsilon$-map onto a graph $G$, then the loop given by $f$ splits as a generator in a free basis in $\pi_1(G)$?
\end{question}   

The answer to this question should be ``yes, of course". However, a second glance at it might change the ``yes, of course" to ``maybe not". Experience shows that one easily gets confused when trying to prove this. As usual, these confusions revolve around a base point change and the affirmative answer to Question \ref{question} starts out being obvious and then it becomes not so obvious and then it becomes just wrong.  

This note is to remark that the answer to Question \ref{question} is no.

\begin{theorem}\label{main theorem}
For each $\varepsilon>0$, there is a surjective $\varepsilon$-map $f\colon \mathbb{S}^1 \rightarrow G$ to a graph $G$ such that $f$ is NOT a generator in a free basis in $\pi_1(G)$. 
\end{theorem} 

This sheds no light on Conjecture \ref{conjecture} but we find it rather surprising. Maybe the reason why this is surprising is because if $\varepsilon$ is small enough, then the loop $f$ is a primitive element in homology. We will leave the proof of this to the readers.

\begin{acknowledgement}
I would like to thank Grigori Avramidi and Max Forester for interesting conversations. I learned the Whitehead algorithm from the second lecture of a three-lecture series by Henry Wilton at the Hausdorff Institute in Fall 2018. I wish the video of his second lecture was available. Finally, I am thankful to the Max Planck Institute for Mathematics for the perfect working environment. 
\end{acknowledgement}

\section{Proof of Theorem \ref{main theorem}}

For each positive integer $k$, take the following word. 
\[w_k = a_1a_2a_1\: a_2a_3a_2\: a_3a_4 a_3\; ... \; a_{2k-1}a_{2k}a_{2k-1}\;\; a_{2k}(\gamma a_1\gamma^{-1}) a_{2k}\; \; \gamma .\]
The word $w_k$ is an element of the free group $F_{2k+1}$ on $(2k+1)$ generators $a_1, a_2,..., a_{2k}$ and $\gamma$. Realize $F_{2k+1}$ as the fundamental group of a wedge of $(2k+1)$ circles that are in one-to-one correspondence with $a_1, a_2,..., a_{2k}, \gamma$. The Whitehead graph\footnote{The definition of the  Whitehead graph is recalled in the Appendix.} of $w$ is given in Figure 1,
\newline 
 
\begin{figure}[h!]
\centering
\includegraphics[scale=0.46]{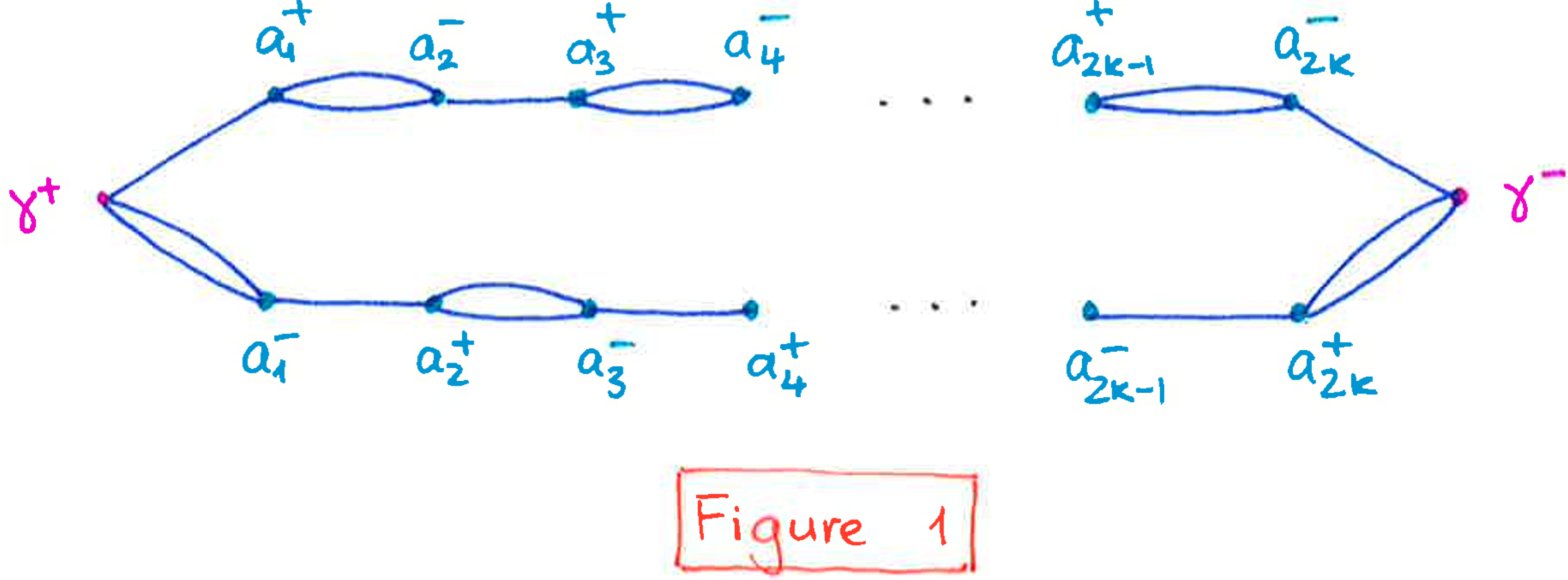}
\end{figure}
\noindent
and does not have a cut vertex\footnote{A vertex $v$ of a graph is a \emph{cut vertex} if the graph consists of two subgraphs with a single vertex in common that is $v$.}. It follows from the Whitehead algorithm (\cite{whiteheadalgorithm}) that $w_k$ is not a generator of $F_{2k+1}$.
We recall the theorem by Whitehead.
\begin{theorem}[Whitehead]
Let $F_n$ be the free group with $n$ generators and let $w$ be a word in the generators. If $w$ is a generator in a free basis of $F_n$, then the Whitehead graph of $w$ has a cut vertex. 
\end{theorem}
\bigskip

Next, we explain how to realize $w_k$ as a surjective  map 
\[f_k\colon(\mathbb{S}^1,x) \rightarrow (G_k, y_1),\] 
for a particular graph $G_k$, such that for a fixed $\varepsilon$, the map $f_k$ is an $\varepsilon$-map if $k$ is large enough. 

The graph $G_k$ is defined as in Figure 2. The base point is $y_1$. The large middle loop, oriented clockwise, is $\gamma$. Each small (embedded) loop $a_i$ based at $y_i$ has length $\pi/(6k)$. To make them based at $y_1$, use the arc going clockwise from $y_1$ to $y_i$ for a change of base point. Each segment $y_iy_{i+1}$ also has length $\pi/(6k)$. 

\begin{figure}[h!]
\centering
\includegraphics[scale=0.5]{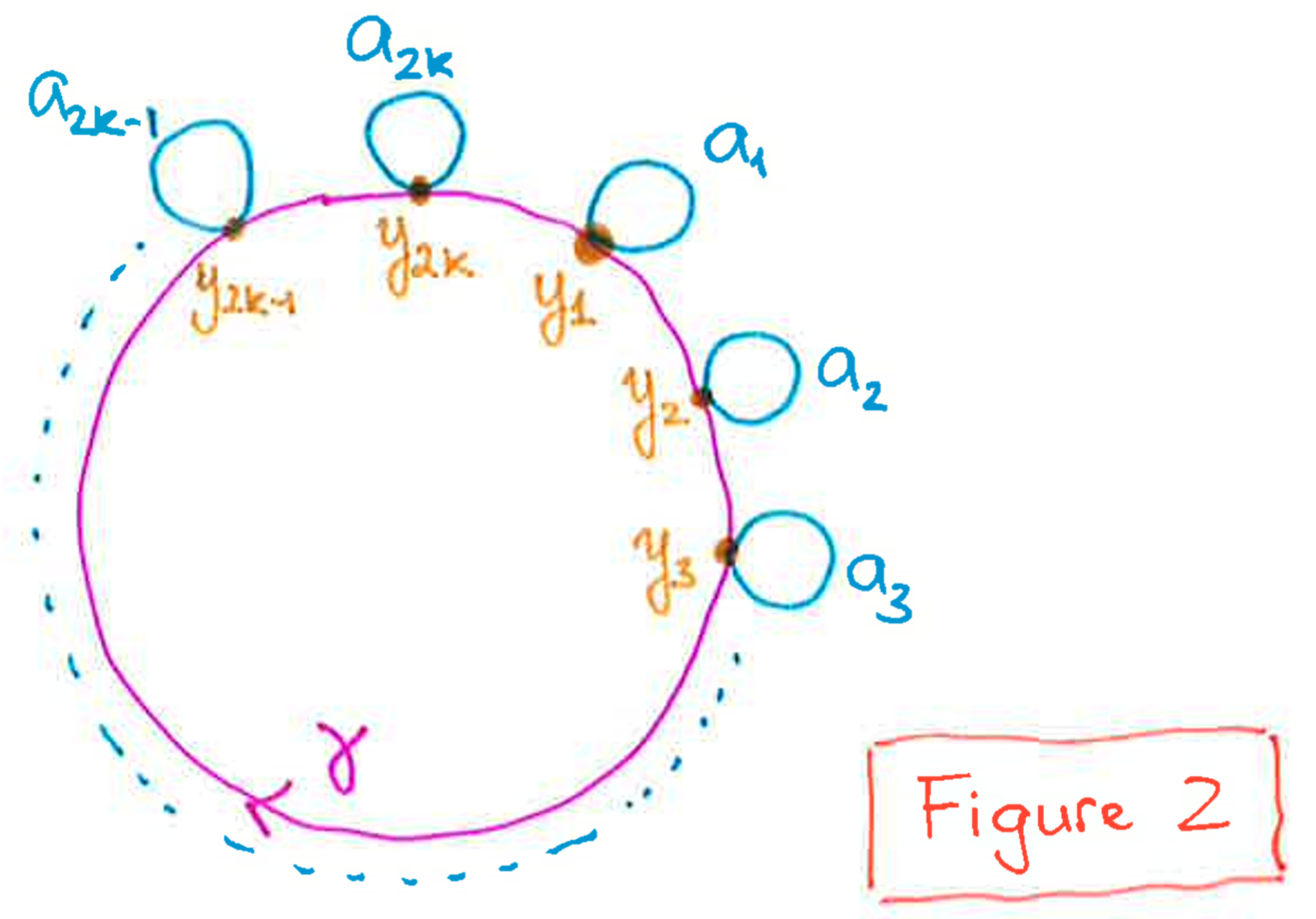}
\end{figure}
The map $f_k$ is constructed as follows. Divide $\mathbb{S}^1$ into $(6\times 2k)$ equal segments. As we go around $\mathbb{S}^1$ starting at the base point $x$, the segments are mapped to $G_k$ in the following order.

\begin{itemize}
\item[.] Loop $a_1$
\item[.] Segment $y_1y_2$
\item[.] Loop $a_2$
\item[.] Segment $y_2y_1$
\item[.] Loop $a_1$
\item[.] Segment $y_1y_2$ (we are done with the segment $a_1a_2a_1$ of $w_k$)
\newline
\item[.] Loop $a_2$
\item[.] Segment $y_2y_3$
\item[.] Loop $a_3$
\item[.] Segment $y_3y_2$
\item[.] Loop $a_2$
\item[.] Segment $y_2y_3$ (we are done with the segment $a_2a_3a_2$ of $w_k$)
\newline 
\item[] ...
\item[] ...

\item[.] Loop $a_{2k}$
\item[.] Segment $y_{2k}y_1$
\item[.] Loop $a_1$
\item[.] Segment $y_1y_{2k}$
\item[.] Loop $a_{2k}$
\item[.] Segment $y_{2k}y_1$ (we are done with $w_k$). 
\end{itemize}

\begin{remark}
The loop we have just defined has cyclic symmetries. However, the element of $\pi_1(G_k) = F_{2k+1}$ does not look like it has cyclic symmetries. This is because when the loop $a_1$ occurs at the end of the journey, we have already basically gone around $\gamma$, so going around $a_1$ then means that $\gamma a_1 \gamma^{-1}$, instead of $a_1$, should be added to the $w$. 
\end{remark}


\section*{Appendix: The Whitehead graph}
Let $F_n$ be the free group on $n$ generators $a_1, a_2, ..., a_n$. Let $w$ be word in the generators. 

The Whitehead graph of $w$ has $2n$ vertices called $a_1^+, a_1^-, a_2^+, a_2^-,..., a_n^+, a_n^-$. To draw the edges of the Whitehead graph, first we put ``$+$" on the right and ``$-$" on the left of each letter that appears in $w$ so that it looks like 
\[ ...\; \; ^{-}a_i^+\quad ^{-}a_j^+\quad ^{-}(a^{-1}_k)^+ \; \; ...\]
For each inverse element $a_k^{-1}$ that occurs in $w$, we change $^{-}(a^{-1}_k)^+$ to $^{+}a_k^-$. Then the chain of symbols we created above becomes
\[ ...\; \; ^{-}a_i^+\quad ^{-}a_j^+\quad ^{+}a_k^- \; \; ...\]
Each pairs of adjacent signs corresponds to an edge in the Whitehead graph. For example, with the above labeling, there is an edge from $a_i^+$ to $a_j^-$ and there is an edge from $a_j^+$ to $a_k^+$. 

One thing one should not forget is that $w$ should be written as a circular word, so there is also an edge from $a_\text{last}^+$ to $a_\text{first}^-$. That is how to draw the Whitehead graph of $w$.

\begin{remark}
The Whitehead graph of $w$ is obtained by first drawing a wedge if $n$ circles intersecting at a point $x_0$ and then drawing the loop given by $w$ in a different color (say, blue) and then cutting out a small neighborhood of $x_0$. The edges are the blue segments of $w$ in this neighborhood and the vertices are where the blue curves exit the neighborhood.  This explains why the word should be treated as circular.
\end{remark}
\bibliography{Reference}
\bibliographystyle{amsplain}

\end{document}